\documentclass[12pt]{amsart}
\include{epsf}
\advance\hoffset by -0.5 truecm
\usepackage{amssymb}
\newtheorem{Theorem}{Theorem}[section]
\newtheorem{Lemma}[Theorem]{Lemma}

\newtheorem{Proposition}[Theorem]{Proposition}

\newtheorem{Definition}[Theorem]{Definition}
\newtheorem{Example}[Theorem]{Example}
\newtheorem{Remark}[Theorem]{Remark}

\def\V{\mbox{Var}}

\def\R\re
\def\V{\bf V}

\def \re{{\mathbb R}}

\def \cp{{\mathbb{CP}}}

\def \0{\lambda_{0}}

\def\comment#1{}

\def\withcomments{
\addtolength{\oddsidemargin}{0.5 in}
\addtolength{\evensidemargin}{0.5 in}
\newcounter{mycommentcounter}
\def\comment##1{\refstepcounter{mycommentcounter}%

\ifhmode%

\unskip%

{\dimen1=\baselineskip \divide\dimen1 by 2 %
\raise\dimen1\llap{\tiny -\themycommentcounter-}}\fi%
\marginpar{\renewcommand{\baselinestretch}{0.8}%
\footnotesize [\themycommentcounter]: \raggedright ##1}}
\date{\framebox{Draft of \today}}}
\withcomments

\begin{document}

\title[On Yamabe constants of Riemannian products]
{On Yamabe constants of Riemannian products}

\author[K.~Akutagawa]{Kazuo Akutagawa}

\thanks{K.~Akutagawa is supported in part by the Grants-in-Aid 
for Scientific Research~(C), 
Japan Society for the Promotion of Science, No.~16540059}
\address{Department of Mathematics, Tokyo University of Science, 
Noda 278-8510, Japan.}
\email{akutagawa${}_{-}$kazuo@ma.noda.tus.ac.jp}
\author[L.~Florit]{Luis A.~Florit}\thanks{L.~A.~Florit is
partially supported by CNPq-FAPERJ, Brazil}
\address{IMPA\\
Est. Dona Castorina, 110\\
22460-320, Rio de Janeiro\\
BRAZIL.}
\email{luis@impa.br}
\author[J.~Petean]{Jimmy Petean}\thanks{J.~Petean is supported by grant 46274-E of CONACYT}
\address{CIMAT  \\
A.P. 402, 36000 \\
Guanajuato. Gto. \\
M\'exico.}
\email{jimmy@cimat.mx}


\date{}

\begin{abstract} 
For a closed Riemannian manifold $(M^m ,g)$ of constant positive scalar curvature 
and any other closed Riemannian manifold $(N^n ,h)$, 
we show that the limit of the Yamabe constants 
of the Riemannian products $(M\times N, g + r h)$ as $r$ goes to infinity
is equal to the Yamabe constant of $(M^m \times \re^n, [g + g_{{}_{\mathbb{E}}}]) $
and is strictly less than the Yamabe invariant of $S^{m+n}$
provided $n\geq 2$.   
We then consider the minimum of the Yamabe functional restricted to functions 
of the second variable and we compute the limit 
in terms of the best constants of the Gagliardo-Nirenberg inequalities. 
\end{abstract}

\maketitle 

\section{Introduction} 
Let $M^m$ be a closed smooth manifold of dimension $m$ 
and denote by $[g]$ the conformal class of a Riemannian metric $g$ on $M$. 
The Yamabe constant $Y(M, [g])$ of $[g]$ 
is the infimum of the normalized total scalar curvature functional 
restricted to $[g]$: 
$$Y(M, [g])= \inf_{h\in [g]}\ 
\frac{\int_M {\bf s}_h d\mu_h}{\textrm{Vol}(M,h)^{\frac{m-2}{m}}} \ \ ,$$ 
\noindent
where ${\bf s}_h$ is the scalar curvature of $h$ and $d\mu_h$ its volume element. 
Of course, the Yamabe constant can be expressed in terms of
functions in the Sobolev space $W^{1,2}(M)$\ \ 
(by writing $h=f^{ \frac{4}{m-2} } g$ if $f \in C_+^{\infty}(M)$):
$$Y(M, [g]) = \inf_{f \in W^{1,2}(M)\atop f \not\equiv 0}Q_g (f) \  
:= \inf_{f\in W^{1,2} (M)\atop f \not\equiv 0}
\frac{\int_M \big(a_m |\nabla f|^2_g  + {\bf s}_gf^2\big) \ 
d\mu_g}{{\| f \| }_{p_m}^2}, $$
\noindent
where $a_m =\frac{4(m-1)}{m-2}$ and $p_m =\frac{2m}{m-2}$.

Functions  realizing the infimum are called {\it Yamabe minimizers} 
and the corresponding metrics are called {\it Yamabe metrics} 
(and have constant scalar curvature). 
They always exist by a fundamental theorem obtained 
in several steps by H. Yamabe, N. Trudinger, T. Aubin and R. Schoen 
\cite{Yamabe, Trudinger, Aubin, Schoen0}. 
Then one defines the Yamabe invariant of $M$, $Y(M)$, 
as the supremum of the Yamabe constants of all conformal classes of metrics 
on $M$ \cite{Schoen, Kobayashi} (cf.~\cite{Kobayashi0}). 
This invariant is always finite since for any conformal class $[g]$, 
$Y(M^m, [g]) \leq Y(S^m, [g_{{}_{S^m}}]) = 
m(m-1) \textrm{Vol}(S^m, g_{{}_{S^m}})^{2/m}$, 
where $g_{{}_{S^m}}$ is the round metric 
on $S^m$ of constant sectional curvature 1. 
We will denote ${\bf Y}_m := Y(S^m) =Y(S^m, [g_{{}_{S^m}}])$. 

The nature of the problem of computing or 
estimating the Yamabe constant of a conformal class, 
and therefore the Yamabe invariant of the manifold, 
depends strongly on whether the constant is 
positive or non-positive. 
If $Y(M, [g]) \leq 0$ then 
$Y(M, [g]) \geq (\inf_M  {\bf s}_g ) \textrm{Vol}(M, g)^{\frac{2}{m}}$, 
as was first pointed out by O. Kobayashi in \cite{Kobayashi}.  
This allows for instance to study the behavior of the Yamabe invariant 
under surgery (\cite{Yun}) 
and so to obtain some understanding of the invariant in the non-positive case; 
\cite{Petean1, Petean2, Rosenberg}. 
Such an estimate does not exist in the positive case. 
In particular there might exist unit volume Riemannian metrics on
$M^m$ of constant scalar curvature greater than ${\bf Y}_m$. 
A typical example of this situation comes from Riemannian products: 
if $(M^m, g)$ and $(N^n, h)$ are unit volume Riemannian manifolds 
of constant scalar curvature and ${\bf s}_g >0$ then, 
for small $\delta > 0$, 
$\delta^n  g + \delta^{-m}  h$ has volume one 
and scalar curvature greater than ${\bf Y}_{m+n}$. 
It is the main purpose of this article 
to study the Yamabe constants of such Riemannian products. 

There is one well understood example in this direction worked out 
by R.~Schoen~\cite{Schoen} and O.~Kobayashi~\cite{Kobayashi0}: 
for any $r > 0$ all solutions of the Yamabe equation on 
$(S^{n-1} \times S^1, g_{{}_{S^{n-1}}} + r  g_{{}_{S^1}})$ 
depend only on the $S^1$-variable 
and one can understand the solutions of the resulting 
ordinary differential equation. 
Following this lead, we will consider Riemannian products 
$\delta^n  g + \delta^{-m} h$ 
on $M \times N$ and look for solutions of the Yamabe equations 
which depend only on the second variable. 

Let 
$$Y_N(M \times N, g + h) := 
\inf_{f \in W^{1,2}(N)\atop f \not\equiv 0}Q_{g+h}(f).$$ 
One can see that the infimum is realized and that if $f$ is such a 
minimizer then $f^{\frac{4}{m+n-2}} \  (g+h)$ has constant scalar curvature 
(we will go over this on Section 2). 
We remark that, 
contrarily to the Yamabe constant $Y(M \times N, [g + h])$, 
this constant $Y_N(M \times N, g + h)$ is not a conformal invariant, 
but merely a scale invariant. 
\vspace{.5cm}

Our first result says in particular that 
the limit of the Yamabe constant of the products above exists: 

\begin{Theorem}\label{3} 
Let $(M^m, g)$ be a closed Riemannian 
$m$-manifold $(m \geq 2)$ of positive scalar curvature 
$($not necessarily constant$)$ 
and $(N^n, h)$ any closed Riemannian $n$-manifold. Then, 
$$\lim_{r \nearrow \infty} Y(M \times N, [g + r h] ) = 
Y(M \times \re^n, [g + g_{{}_{\mathbb{E}}}]) 
:= \inf_{f\in C_c^{\infty}(M \times \re^n) \atop f \not\equiv 0}
Q_{g + g_{{}_{\mathbb{E}}}}(f) > 0. $$ 
\noindent 
And 
$$\lim_{r \nearrow \infty} Y_N(M \times N, g + r h ) = 
Y_{\re^n}(M \times \re^n, g + g_{{}_{\mathbb{E}}}) 
:= \inf_{f\in C_c^{\infty}(\re^n) \atop f \not\equiv 0}
Q_{g + g_{{}_{\mathbb{E}}}}(f) > 0. $$ 
Here, $g_{{}_{\mathbb{E}}}$ stands for the Euclidean metric on $\re^n$. 
\end{Theorem} 

\begin{Remark} 
{\rm The fact that $Y(S^{n-1} \times S^1) = Y(S^n)$ (for $n \geq 3$) was first proved by 
O.~Kobayashi~\cite{Kobayashi0} and R.~Schoen~\cite{Schoen}, 
by analysis of the behavior of constant scalar curvature metrics 
on $(\mathbb{R}^n - \{0\}, g_{{}_{\mathbb{E}}})$. 
Another proof was given by also O.~Kobayashi~\cite{Kobayashi}, 
by using an argument in the proof of 
the celebrated Kobayashi's inequality~\cite[Theorem~2]{Kobayashi}. 
The above theorem gives the third proof since 
$$Y(S^{n-1} \times S^1) \geq 
Y(S^{n-1} \times \mathbb{R}^1, [g_{{}_{S^{n-1}}} + g_{{}_{\mathbb{E}}}]) 
= Y(S^n, [g_{{}_{S^n}}])\ \bigl(= Y(S^n)\bigr).$$ }
\end{Remark} 

On the constant $Y(M \times \re^n, [g + g_{{}_{\mathbb{E}}}])$, 
we also obtain:  

\begin{Theorem}\label{4} 
Let $(M^m, g)$ be a closed Riemannian $m$-manifold $(m \geq 2)$ 
of positive scalar curvature $($not necessarily constant$)$. 
Assume that $n \geq 2$. Then, 
$$ 
Y(M \times \re^n, [g + g_{{}_{\mathbb{E}}}]) < {\bf Y}_{m + n}. 
$$ 
\end{Theorem} 

Recall that the best $n$-dimensional Sobolev constant 
is the smallest positive number $\sigma_n$ such that 
for any smooth compactly supported function $f$ on 
$\re^n$, ${\| f \|}_{p_n }^2 \leq \sigma_n {\| \nabla f \|}_2^2 $. 
Due to the conformal invariance of the Yamabe constant, 
one can use the stereographic projection to translate the Yamabe functional 
from the round sphere to the Euclidean space to obtain:
$$\sigma_n =\frac{a_n}{{\bf Y}_n}.$$ 

In a similar fashion we will see that when studying the limits above 
a fundamental role is played by the best constant in the Gagliardo-Nirenberg inequalities: 
namely, we will call $\sigma_{m,n}$ the smallest positive number 
such that for any $f\in W^{1,2} ({\re}^n)$,
$$ {\| f \|}_{p_{m+n}}^2      \leq    \sigma_{m,n} \  
{\| \nabla f \|}_2^{\frac{2n}{m+n}} {\| f \|}_2^{\frac{2m}{m+n}}.$$ 
\noindent 
Or what is equivalent: 
$$\sigma_{m,n} = {\left( \inf_{f\in W^{1,2}(\re^n) \atop f \not\equiv 0} 
\frac{{\| \nabla f \|}_2^{\frac{2n}{m+n}} 
{\| f \|}_2^{\frac{2m}{m+n}}}{{\| f \|}_{p_{m+n}}^2   } \right) }^{-1}. $$ 
\noindent 
The constant $\sigma_{m,n}$ is of classical interest in the study of 
partial differential equations and has been computed numerically, 
although it is not known any explicit expresion for the constant 
or for the minimizing function (which is known to be radial 
and decreasing); \cite{g1, g2, n, k, Levine, Weinstein}. 
\vspace{0.2cm} 

We will prove: 

\begin{Theorem}\label{1} 
Let $(M^m,g)$ be a closed smooth unit volume 
Riemannian manifold of constant positive scalar curvature ${\bf s}_g$. 
Then 
$$Y_{\re^n}(M \times \re^n, g + g_{{}_{\mathbb{E}}}) 
\ = \ \frac{C(m,n) \ {\bf s}_g^{\frac{m}{m+n}} }{ \sigma_{m,n}}, $$ 
where $C(m,n) = (a_{m+n} )^{\frac{n}{m+n}} (m+n) n^{\frac{-n}{m+n}}m^{\frac{-m}{m+n}}$. 
\end{Theorem} 
\vspace{0.1cm}

It is clear that 
$$Y(M \times \re^n, [g + g_{{}_{\mathbb{E}}}]) 
\leq Y_{\re^n}(M \times \re^n, g + g_{{}_{\mathbb{E}}})$$ 
and it seems that equality should hold under certain hypothesis. 
It certainly cannot always be the case since 
$Y_{\re^n}(M \times \re^n, g + g_{{}_{\mathbb{E}}}) > {\bf Y}_{m+n}$ 
if ${\bf s}_g$ is big enough. 
\vspace{.3cm}
\noindent

{\bf Question :} 
Is it true that 
$Y_{\re^n}(M \times \re^n, g + g_{{}_{\mathbb{E}}}) 
= Y(M \times \re^n, [g + g_{{}_{\mathbb{E}}}])$ 
if $g$ is a Yamabe metric ? 
\vspace{.3cm}
\noindent

As we mentioned before, 
the constants $\sigma_{m,n}$ can be explicitly computed numerically. 
Using these computations, 
we apply Theorem \ref{1} to the case 
when $(M,g)$ and $(N,h)$ are round spheres. 
These are particularly interesting cases 
because of Schoen and Kobayashi's argument mentioned above 
and because $S^n \times S^m$ is obtained
by performing surgery on $S^{m+n}$, 
and therefore if the surgery theorem in \cite{Yun}
were true in the positive case 
we should have $Y(S^m \times S^n ) = {\bf Y}_{m+n}$.
Set 
$$Y^{\infty}_{S^n}(S^m \times S^n ) 
:= \lim_{r \nearrow \infty} Y_{S^n} (S^m \times S^n, g_{{}_{S^m}} + r g_{{}_{S^n}}) 
= Y_{\re^n}(S^m \times \re^n, g_{{}_{S^m}} + g_{{}_{\mathbb{E}}}).$$ 
We give the corresponding values for all $m,n\geq 2$ with
$m+n\leq 9$. 
$$ 
\begin{array}{lllll}\label{lista} 
m & n & \ \ \sigma_{m,n}^{-1} & Y^{\infty}_{S^n}(S^m\!\times\!S^n) & \ \ {\bf Y}_{m+n}\\ 
\hline 
\vspace{.2cm} 

2 & 2 & 2.41877 &\ \ \ 59.40481 & 61.56239\\ 
2 & 3 & 3.87947 &\ \ \ 75.39687 & 78.99686\\ 
\vspace{.2cm} 

3 & 2 & 2.11360 &\ \ \ 78.18644 & 78.99686\\ 
2 & 4 & 5.66408 &\ \ \ 91.68339 & 96.29728\\ 
3 & 3 & 3.19925 &\ \ \ 94.71444 & 96.29728\\ 
\vspace{.2cm} 

4 & 2 & 1.90282 &\ \ \ 95.87367 & 96.29728\\ 
2 & 5 & 7.71937 &\ \ \ 108.1625 & 113.5272\\ 
3 & 4 & 4.53960 &\ \ \ 111.2934 & 113.5272\\ 
4 & 3 & 2.75810 &\ \ \ 112.6214 & 113.5272\\ 
\vspace{.2cm} 

5 & 2 & 1.75469 &\ \ \ 113.2670 & 113.5272\\ 
2 & 6 & 10.0021 &\ \ \ 124.7747 & 130.7157\\ 
3 & 5 & 6.10843 &\ \ \ 127.9414 & 130.7157\\ 
4 & 4 & 3.81586 &\ \ \ 129.3551 & 130.7157\\ 
5 & 3 & 2.45567 &\ \ \ 130.1272 & 130.7157\\ 
\vspace{.2cm} 

6 & 2 & 1.64650 &\ \ \ 130.5398 & 130.7157\\ 
2 & 7 & 12.4764 &\ \ \ 141.4740 & 147.8778\\ 
3 & 6 & 7.88171 &\ \ \ 144.6521 & 147.8778\\ 
4 & 5 & 5.06274 &\ \ \ 146.1089 & 147.8778\\ 
5 & 4 & 3.32083 &\ \ \ 146.9519 & 147.8778\\ 
6 & 3 & 2.23778 &\ \ \ 147.4615 & 147.8778\\ 
7 & 2 & 1.56455 &\ \ \ 147.7507 & 147.8778 
\end{array}
$$ 
\vspace{.3cm}

\noindent It should be the case that 
$Y_{S^n}^{\infty}(S^m \times S^n ) < {\bf Y}_{n+m}$ 
for all values $m,n \geq 2$. 
Hence, this gives a proof of Theorem~\ref{4} 
for $(M^m, g) = (S^m, g_{{}_{S^m}})$, 
namely that of 
$$Y(S^m \times \re^n, [g_{{}_{S^m}} + g_{{}_{\mathbb{E}}}]) < {\bf Y}_{m + n}$$ 
when $m, n \geq 2, m + n \leq 9$. 
Moreover, the above method gives a numerical estimate from above for the constant 
$Y(S^m \times \re^n, [g_{{}_{S^m}} + g_{{}_{\mathbb{E}}}])$. 

Note also that in the $4$--dimensional case the Yamabe invariant of 
$\cp^2$ is realized by the conformal class of the Fubini-Study 
metric $g_{\textrm{FS}}$ (\cite{LeBrun, Gursky}), 
giving 
$$Y(\cp^2) = Y(\cp, [g_{\textrm{FS}}]) = 12 \sqrt{2} \pi = 53.31459\cdots$$ 
\noindent 
and since Einstein metrics are always Yamabe metrics we have that 
$$Y(S^2 \times S^2, [g_{{}_{S^2}} + g_{{}_{S^2}}])=16 \pi = 50.26548\cdots.$$ 

In the next section, 
we will review some known results on Yamabe constants, 
point out a few observations and fix some notation. 
In Section 3, we will recall Schoen and Kobayashi's discussion 
of the solutions of the Yamabe equation on 
$(S^{n-1} \times S^1, g_{{}_{S^{n-1}}} \times r g_{{}_{S^1}})$. 
We will prove Theorem~\ref{3} and Theorem~\ref{4} in Section 4 
and Theorem~\ref{1} in Section 5. 
Finally, a procedure to compute 
numerically these Yamabe constants is given in the last section. 
\vspace{1cm}

{\it Acknowledgements:} 
The authors would like to thank Professors Manuel del Pino and Jean Dolbeault 
for valuable observations on the Gagliardo-Niremberg inequalities 
and Fernando Coda for numerous conversations on the subject. 
The first author would like to express his gratitude to 
Professors Boris Botvinnik and Osamu Kobayashi for many helpful discussions. 
The third author would like to express his gratitude to IMPA 
where this work was carried on with the partial support of CAPES-Brazil.

\section{Preliminaries} 
Let $(X^k ,g)$ be a closed smooth $k$-dimensional Riemannian manifold. 
Recall that  ${\bf s}_g$ is 
the scalar curvature of $g$, $d\mu_g$ its volume element and 
$$p=p_k = \frac{2k}{k - 2} \ \ \ \ \ \ \mbox{and}\ \ \ \ \ \ 
a = a_k = \frac{4(k - 1)}{k - 2}.$$ 
Consider the Sobolev space $W^{1,2}(X)$ and 
the {\it Yamabe functional} defined by 
$$f \hspace{.5cm} \mapsto \hspace{.5cm} 
Q_g(f) := 
\frac{\int_X \big(a_k|\nabla f|_g^2 d\mu_g \ + \ {\bf s}_gf^2\big) d\mu_g}
{{\| f \|}_p^2 }.$$ 
We say that $f$ is a {\it Yamabe minimizer} 
(for $g$) if it realizes the minimum of the Yamabe functional. 
In this case $f^{4/(k-2)} g $ has constant scalar curvature and the Yamabe constant 
of the conformal class of $g$ is then equal to $Q_g (f)$. 

In this paper we consider a unit volume Riemannian product 
$(M^m \times N^n, g + h)$. 
We assume that the scalar curvature of both $g$ and
$h$ are constant and try to understand the Yamabe constant 
of the conformal class of the product metric. 
This is of no interest if ${\bf s}_g + {\bf s}_h$ is negative, 
since in this case we have uniqueness of the Yamabe metric. 
The situation we want to address is when 
${\bf s}_g + {\bf s}_h$ is positive and bigger than ${\bf Y}_{m+n}$, 
the Yamabe invariant of the round sphere $S^{m+n}$. 
In this case there must exist a non-constant Yamabe function, 
and so another metric of constant scalar curvature  in the same conformal class. 

To compute the Yamabe constant is a very difficult problem 
and so it is to understand the Yamabe minimizer. 
We will then restrict ourselves to functions 
which depend only on one of the variables, that is, 
positive smooth functions $f : N^n \rightarrow \re$ of one 
of the factors in the Riemannian product $M^m \times N^n$. 
The scalar curvature of $f^{\frac{4}{m+n-2}}(g+h)$ is given by 
$${\bf s}_{g+h} = f^{1-p_{m+n}} \ (-\ a_{m+n}{\Delta}_hf\ + \ {\bf s}_{g+h}f).$$ 

We then introduce the following definition:

\begin{Definition} 
For a given Riemannian product $(M \times N, g + h)$ of 
constant scalar curvature, 
we will call the $N$--{\rm Yamabe constant} to 
the infimum of the $(g+h)$--Yamabe functional restricted to $W^{1,2}(N)$. 
We will denote this constant by $Y_N(M \times N, g + h)$. 
\end{Definition} 

To study critical points of the $(g+h)$--Yamabe functional restricted to $W^{1,2}(N)$, 
let $\varphi ,f: N\rightarrow \re$ be smooth functions. 
A well-known computation gives that 
$$\frac{d\big(Q(f + t\varphi)\big)}{dt}(0) = 
\frac{2 \textrm{Vol}(M,g)}{{\| f \| }_p^2 } \int_{N}\ \Big(-\ a_{m+n}
{\Delta}_h f + {\bf s}_{g+h} f - {\| f \|}_p^{-p+2}Q(f) f^{p-1} \Big)
\ \varphi \ d\mu_h,$$ 
where $p=p_{m+n}$ and $Q(\cdot) = Q_{g+h}(\cdot)$. 
Therefore the critical points of the Yamabe functional restricted to 
$W^{1,2}(N)$ are precisely the functions $f$ such that 
the conformal metric $f^{4/(m+n-2)}(g+h)$ has constant scalar curvature 
${\bf \tilde{s}} = {\| f \|}_p^{-p+2}Q(f)$. 
The next point is that the infimum of the Yamabe functional 
restricted to $W^{1,2}(N)$ is always achieved. 
This is a simple fact, 
it is essentially the subcritical case of the Yamabe problem, 
but we sketch its proof 
since we have not seen it in the literature. 
\vspace{0.1cm}

\begin{Proposition}\label{easy} 
Let $(M^m ,g)$ and $(N^n ,h)$ be closed Riemannian manifolds of constant scalar curvature. 
Then, there exists a positive smooth function $f$ on $N^n$ 
which minimizes the Yamabe quotient among all functions on $N^n$. 
\end{Proposition} 

\begin{proof} 
Without loss of generality, 
we may assume that $\textrm{Vol}(M, g) = 1$. 
Let $\{u_i\}$ be a sequence of non-negative functions on $N$ 
such that $Q_{g + h}(u_i) \rightarrow Y_N(M \times N, g + h)$. 
We can assume that ${\| u_i \|}_{p_k} =1$, 
where $k = m + n$. 
Then, 
$$ {\| u_i \|}^2_{2,1} = \int_N \big(|\nabla u_i |^2_h + u_i^2\big) \ d\mu_h $$ 
$$\ \ \ \  \ \ \ \ \ \ \  \ \ \ \ \ \ \ \ \ \ \ \ \ 
=\frac{1}{a_k} Q_{g + h}(u_i ) + 
\int_N \big(1 - \frac{{\bf s}_{g+h}}{a_k}\big)u_i^2 \ d\mu_h $$ 
$$\ \ \   \leq \frac{Y_N(M \times N, g + h) + 1}{a_k} + K{\| u_i \|}^2_2$$ 
for some $K > 0$, 
that is bounded independently of $i$ 
since ${\| u_i \|}^2_2 \leq {\| u_i \|}^2_{p_k}
\textrm{Vol}(N, h)^{2/k} = \textrm{Vol}(N, h)^{2/k}$ 
by H\"{o}lder's inequality. 
Since 
$$\frac{1}{p_k} > \frac{1}{2} - \frac{1}{n},$$
by the Rellich-Kondrakov Theorem 
the inclusion $W^{1,2}(N) \subset L^{p_k}(N)$ is a compact operator. 
Therefore, there exists a subsequence of the $\{u_i\}$ 
which converges weakly in $W^{1,2}(N)$ and strongly in $L^{p_k}(N)$ 
to a function $u \in W^{1,2}(N)$ with ${\| u \|}_{p_k} = 1$. 

Now, by the weak convergence in $W^{1,2}(N)$, we have 
$${\big|\big| |\nabla u| \big|\big|}_2^2 = \lim_{i \rightarrow \infty}
\int_N \langle \nabla u,\nabla u_i \rangle \ d\mu_h, $$ 
and therefore 
$${\big|\big| |\nabla u| \big|\big|}_2^2 \leq \limsup_{i \rightarrow \infty} 
{\big|\big| |\nabla u_i | \big|\big|}_2^2. $$ 
And since by H\"{o}lder's inequality 
$$\int_N {\bf s}_{g + h}u^2 d\mu_h = 
\lim_{i \rightarrow \infty} \int_N {\bf s}_{g + h}u_i^2 d\mu_h, $$ 
we have that $Q_{g + h}(u) \leq \lim_{i \rightarrow \infty} Q_{g + h}(u_i)$, 
and hence $Q_{g + h}(u) = Y_N(M \times N, g + h)$. 
It then follows from elliptic regularity 
that $u$ has to be strictly positive and smooth. 
\end{proof} 

\begin{Remark} 
{\rm Note that, 
for a given Riemannian product of constant scalar curvature, 
we have three associated Yamabe constants 
each producing a constant scalar curvature metric. 
The three are equal if the original product is a Yamabe metric.} 
\end{Remark}

\section{Reviewing the circle}

Schoen~\cite{Schoen} (cf.~Kobayashi~\cite{Kobayashi0}) 
gave a fairly complete study of the solutions of the Yamabe equation for the manifolds 
$(S^{n-1} \times S^1, g_{{}_{S^{n-1}}} + r g_{{}_{S^1}})$, 
where $n \geq 2$ and $r$ is a positive constant. 
He points out that due to the conformal invariance 
and a theorem of Caffarelli--Gidas--Spruck \cite{Caffarelli} 
all solutions are functions of $S^1$. 
Moreover, he writes down the Yamabe equation for a function of $S^1$. 
Moving to the Riemannian universal covering 
$(S^{n-1} \times \re, g_{{}_{S^{n-1}}} + dt^2)$, 
one has to deal with the ordinary differential equation: 

$$\frac{d^2 u}{dt^2} - \frac{1}{4}(n-2)^2 u +\frac{n(n-2)}{4}
u^{(n+2)/(n-2)} = 0,$$ 

\noindent
and look for solutions which are $2\pi r$--periodic. 
Note that exactly the same equation shows up 
if we consider a Riemannian product 
$(M \times S^1, g + r g_{{}_{S^1}})$, 
where $M$ is $(n-1)$-dimensional and the scalar
curvature of $g$ is $(n-1)(n-2)$. 
In this way, 
one can understand all constant scalar curvature metrics 
which are conformal by a function of $S^1$ to 
$(M \times S^1, g + r g_{{}_{S^1}})$; 
the solutions are the same as those for $S^{n-1}$ 
discussed in \cite{Schoen} and \cite{Kobayashi0}. 
So for $r$ close to 1 there is going to be only one solution, 
and as $r$ increases the number of solutions will increase. 
If $\textrm{Vol}(M,g) = V_{n-1} := \textrm{Vol}(S^{n-1}, g_{{}_{S^{n-1}}})$, 
then the $S^1$-Yamabe constant of the product 
will be less than ${\bf Y}_n$ for all $r$ 
and will approach ${\bf Y}_n$ as $r$ goes to infinity. 

\begin{Example} 
{\rm If $\textrm{Vol}(M,g) > V_{n-1}$ in the discussion above, 
the number of solutions will still be the same, 
but as $r$ becomes big the $S^1$--Yamabe constants of the product 
will be bigger than ${\bf Y}_n$. 
In particular, the $S^1$--Yamabe constant will be greater than the Yamabe constant.} 
\end{Example}

\section{Proofs of Theorem~\ref{3} and Theorem~\ref{4}}

\begin{proof}[Proof of Theorem~\ref{3}] 
To simplify the notation, 
we set 
$g_r := g + r h$ on $M^m \times N^n$ and 
$\bar{g} := g + g_{{}_{\mathbb{E}}}$ on $M^m \times \re^n$. 
We may also assume that $\textrm{Vol}(M, g) = 1$. 

First, we show the following 
$$Y(M \times {\mathbb R}^n, [\bar{g}]) > 0.$$ 
Note that $(M \times \re^n, \bar{g})$ is a complete Riemannian manifold 
with strictly positive injective radius 
and bounded sectional curvature. 
Under this conditions, 
the Sobolev embedding 
$W^{1,2}(M \times \re^n) \hookrightarrow L^p(M \times \re^n)$ holds 
(cf.~\cite[Theorem~2.21]{Aubin-Book}), 
that is, there exists constant $K_1 > 0$ such that 
$$\|f\|^2_p \leq K_1\|f\|^2_{2,1}\quad 
\textrm{for}\ f \in W^{1,2}(M \times {\mathbb R}^n),$$ 
where $p = p_{m + n} := \frac{2(m + n)}{m + n - 2}$. 
This and the positivity of the scalar curvature of $(M, g)$ imply that 
$$\Bigl(\int_{M \times {\mathbb R}^n}|f|^pd\mu_{\bar{g}}\Bigr)^{2/p} 
\leq  \frac{K_1}{\alpha}\int_{M \times {\mathbb R}^n}
\bigl(a_{m+n}|\nabla f|^2_{\bar{g}} + {\bf s}_gf^2\bigr)d\mu_{\bar{g}}$$ 
for $f \in W^{1,2}(M \times {\mathbb R}^n),$ where 
$\alpha := \min \bigl\{a_{m+n},\ \underset{M}{\min}\ {\bf s}_g\bigr\} > 0,$ 
and hence 
$$Y(M \times {\mathbb R}^n, [\bar{g}]) \geq \frac{\alpha}{K_1} > 0.$$ 
We also have 
$$Y_{\re^n}(M \times {\mathbb R}^n, \bar{g}) \geq 
Y(M \times {\mathbb R}^n, [\bar{g}]) > 0.$$ 

Second, we prove the following 
\begin{equation}\label{I} 
\liminf_{r \nearrow \infty}Y(M \times N, [g_r]) 
\geq Y(M \times {\mathbb R}^n, [\bar{g}]) 
\end{equation} 
\noindent 
and 
\begin{equation}\label{II} 
\liminf_{r \nearrow \infty}Y_N(M \times N, g_r ) 
\geq Y_{{\mathbb R}^n}(M \times {\mathbb R}^n, \bar{g}) \ . 
\end{equation} 
Pick any $\varepsilon > 0$. 
There exist a small constant $\delta > 0$ and 
finite points $\{q_1, \cdots, q_{\ell}\} \subset N$ such that 
$$\bullet\quad \{U_k := 
\textrm{exp}^h_{q_k}\bigl({\bf B}_{\delta}^h(0)\bigr)\}_{k = 1}^{\ell}\quad 
\textrm{is an open covering of}\ \ N 
\qquad \qquad \qquad \qquad \qquad \qquad \qquad \qquad \qquad \qquad$$ 
and that, on each $U_k$ with respect to 
$h$-normal coordinates $x = (x^1, \cdots, x^n)$ at $q_k$, 
$$\bullet\quad (1 + \varepsilon)^{-1}\delta_{ij}dx^idx^j 
\leq h_{ij}dx^idx^j \leq (1 + \varepsilon)\delta_{ij}dx^idx^j, 
\qquad \qquad \qquad \qquad \qquad \qquad \qquad \qquad \qquad \qquad $$ 
$$\bullet\quad (1 + \varepsilon)^{-1}dx
\leq d\mu_h \leq (1 + \varepsilon)dx. 
\qquad \qquad \qquad \qquad \qquad \qquad \qquad \qquad \qquad \qquad 
\qquad \qquad \qquad \qquad \qquad $$ 
Here, 
$\textrm{exp}^h_{q_k} : {\bf B}_{\delta}^h(0) 
:= \{v \in T_{q_k}N\ \big|\ |v|_h < \delta\} 
\rightarrow U_k$ and $dx := dx^1\wedge\cdots \wedge dx^n$ 
denote respectively the $h$-exponential map at $q_k$ 
and the Euclidean volume form. 
Then note that, 
for any $r > 1$, on each $U_i$ 
with respect to $(r^2 h)$-normal coordinates $y = (y^1, \cdots, y^n)$ at $q_k$, 
$$(1 + \varepsilon)^{-1}\delta_{ij}dy^idy^j 
\leq (r^2 h)_{ij}dy^idy^j \leq (1 + \varepsilon)\delta_{ij}dy^idy^j, $$ 
$$(1 + \varepsilon)^{-1}dy
\leq d\mu_{r^2 h} \leq (1 + \varepsilon)dy. $$ 
We also note that there exists a constant $K_2 > 0$ such that 
$$\textrm{diam}(U_k, r^2 h) \geq K_2 \ r\quad (k = 1, \cdots, \ell)$$ 
for any $r \geq 1$. 
Let $\{\eta_k = \chi_k^2\}_{k = 1}^{\ell}$ be 
a partition of unity subordinate to the covering $\{U_k\}_{k = 1}^{\ell}$ 
and  $K_3> 0$ a positive constant independent of $r \geq 1$
such that 
$$|\nabla \chi_k|_h \leq K_3\qquad (k = 1, \cdots, \ell),$$ 
and hence 
$$|\nabla \chi_k|_{r^2 h} \leq \frac{K_3}{r}\qquad (k = 1, \cdots, \ell).$$ 

With the above understandings, 
for any $r > 1$ and $\varphi \in C^{\infty}(M \times N)$, 
we estimate the $L^p$-norm of $\varphi$ with respect to $g_r$ as follows: 

\begin{align*} 
||\varphi||^2_p 
& = ||\varphi^2||_{p/2} 
= ||\Sigma_k\chi_k^2\varphi^2||_{p/2} \\ 
& \leq \Sigma_k 
\Bigl(\int_{M \times U_k}|\chi_k \varphi|^{p}d\mu_{g_r})\Bigr)^{2/p} \\ 
& \leq (1 + \varepsilon)^{2/p}\Sigma_k 
\Bigl(\int_{M \times U_k}|\chi_k \varphi|^{p}d\mu_{\bar{g}}\Bigr)^{2/p}. 
\end{align*} 

Here, we identify $U_k = \textrm{exp}^{r^2 h}_{q_k}
\bigl({\bf B}^{r^2 h}_{r\delta}(0)\bigr)\ (\subset N)$ 
with $\mathbb{B}_{r\delta}(0) := \{y \in \mathbb{R}^n\ \big|\ |y| < r\delta\}$ 
via the composition of the inverse $\bigl(\textrm{exp}^{r^2 h}_{q_k}\bigr)^{-1}$ 
and the identification ${\bf B}^{r^2 h}_{r\delta}(0) \cong \mathbb{B}_{r\delta}(0)$. 

Set $Y_0 := Y(M \times \re^n, [\bar{g}]) $ and 
$Y_0^0 := Y_{\re^n}(M \times \re^n, \bar{g} ) $. 
We also note that, on each $M \times U_k\ (\subset M \times \mathbb{R}^n)$, 

\begin{align*} 
& \Bigl(\int_{M \times U_k}|\chi_k \varphi|^p d\mu_{\bar{g}}\Bigr)^{2/p} 
\leq \frac{1}{Y_0}
\Bigl(a_{m+n}\int_{M \times U_k}|\nabla (\chi_k \varphi)|_{\bar{g}}^2d\mu_{\bar{g}} 
+ \int_{M \times U_k}{\bf s}_g \chi_{\ell}^2\varphi^2d\mu_{\bar{g}} \Bigr) \\ 
& \qquad \qquad \leq \frac{(1 + \varepsilon)^2}{Y_0}
\Bigl(a_{m+n}\int_{M \times U_k}|\nabla (\chi_k \varphi)|_{g_r}^2d\mu_{g_r} 
+ \int_{M \times U_k}{\bf s}_g \chi_k^2\varphi^2d\mu_{g_r} \Bigr) \\ 
& \qquad \qquad \leq \frac{(1 + \varepsilon)^3}{Y_0}
\Biggl(a_{m+n}\int_{M \times U_k}\chi_k^2|\nabla \varphi|_{g_r}^2d\mu_{g_r} 
+ \int_{M \times U_k}{\bf s}_g \chi_k^2\varphi^2d\mu_{g_r} \\ 
& \qquad \qquad \qquad \qquad \qquad \qquad \qquad \qquad \qquad \qquad \ 
+ \frac{a_{m+n}K_3^2(1 + \varepsilon^{-1})}{r^2}
\int_{M \times U_k}\varphi^2d\mu_{g_r} \Biggr). 
\end{align*} 

\noindent 
Here, it is important to note for the proof of (\ref{II}) that, 
if $\varphi \in C^{\infty}(N)$, then we can replace $Y_0$ by $Y_0^0$. 

Hence, we have 

\begin{align*} 
& \quad \ ||\varphi||^2_p \leq \\ 
&  \frac{(1 + \varepsilon)^{3 + (2/p)}}{Y_0}
\Biggl(\int_{M \times N}\Bigl(a_{m+n}|\nabla \varphi|_{g_r}^2 
+ {\bf s}_g \varphi^2\Bigr)d\mu_{g_r} 
+ \frac{\ell a_{m+n}K_3^2(1 + \varepsilon^{-1})}{r^2}
\int_{M \times N}\varphi^2d\mu_{g_r} \Biggr). 
\end{align*}

 From the positivity of the scalar curvature ${\bf s}_g$, 
there exists a large constant 
$r_0 = r_0(\varepsilon, \underset{M}{\textrm{min}}~{\bf s}_g , (N, h), m+n) > 1$ such that 
$$\frac{\ell a_{m+n}K_3^2(1 + \varepsilon^{-1})}{r_0^2} 
\leq \bigl(\underset{M}{\textrm{min}}~ {\bf s}_g \bigr) \varepsilon.$$ 
Therefore, we obtain 
$$||\varphi||^2_p 
\leq \frac{(1 + \varepsilon)^{4 + (2/p)}}{Y_0} 
\int_{M \times N}\Bigl(a_{m+n}|\nabla \varphi|_{g_r}^2 
+ {\bf s}_g \varphi^2\Bigr)d\mu_{g_r} $$ 
for any $r \geq r_0$ and $\varphi \in C^{\infty}(M \times N)$. 
Again, if  $\varphi \in C^{\infty}(N)$, we can replace $Y_0$ by $Y_0^0$. 
Then, this implies that 
$$Y(M \times N, [g_r]) \geq 
\frac{Y_0}{(1 + \varepsilon)^{4 + (2/p)}} \ ,\qquad 
Y_N (M \times N, g_r ) \geq 
\frac{Y_0^0}{(1 + \varepsilon)^{4 + (2/p)}}$$ 
for any $r \geq r_0$. And since $\varepsilon > 0$ is arbitrary, 
$$\liminf_{r \nearrow \infty}Y(M \times N, [g_r]) \geq Y_0 $$ 

\noindent 
and 
$$\liminf_{r \nearrow \infty}Y_N(M \times N, g_r ) \geq Y_0^0 .$$ 

Finally, we prove 

\begin{equation} 
\limsup_{r \nearrow \infty}Y(M \times N, [g_r]) 
\leq Y_0
\end{equation} 

\noindent
and 

\begin{equation}\label{III} 
\limsup_{r \nearrow \infty}Y_N (M \times N, g_r )  \leq Y_0^0 . 
\end{equation} 

Note that 

\begin{equation} 
\lim_{\rho \nearrow \infty}Y(M \times \mathbb{B}_{\rho}(0), [\bar{g}]) 
= Y_0, 
\end{equation}  

\noindent 
and 

\begin{equation} 
\lim_{\rho \nearrow \infty}
Y_{\mathbb{B}_{\rho}(0)}(M \times \mathbb{B}_{\rho}(0), [\bar{g}]) 
= Y_0^0 . 
\end{equation}  

\noindent 
Take any small $\varepsilon > 0$ and large $\rho > 0$. 
Fix a point $q \in N$ and set an $h$-normal open neighborhood 
$U := \textrm{exp}^h_q\bigl({\bf B}_{\delta}^h(0)\bigr)$ of $q$ 
for $\delta > 0$. 
Here, we choose $\delta = \delta(\varepsilon, (N, h)) > 0$ sufficiently small 
satisfying the same conditions on $U$ as those in the preceding argument. 
Let $r_1 > 0$ be a positive constant such that $r_1\delta \geq \rho$. 
For each $r \geq r_1$, we also use 
the $(r^2 h)$-normal coordinates $y = (y^1, \cdots, y^n)$ at $q$ on $U$ 
and the identification 
$$\mathbb{B}_{r\delta}(0)\ (\subset \mathbb{R}^n) 
\overset{\cong}{\longleftrightarrow} 
U = \textrm{exp}^{r^2 h}_q\bigl({\bf B}_{r\delta}^{r^2 h}(0)\bigr)\ (\subset N).$$  

With the above understandings, 
for any $r \geq r_1$ and 
$f \in C^{\infty}_c(M \times \mathbb{B}_{\rho}(0))\\ 
\Bigl(\subset C^{\infty}_c(M \times \mathbb{B}_{r\delta}(0)) 
\cong C^{\infty}_c(M \times U)\Bigr)$, we obtain 

\begin{align*} 
||f||^2_p
& = \Bigl(\int_{M \times \mathbb{B}_{\rho}(0)}|f|^p d\mu_{\bar{g}}\Bigr)^{2/p} 
\leq(1 + \varepsilon)^{2/p}
\Bigl(\int_{M \times U}|f|^p d\mu_{g_r} \Bigr)^{2/p} \\ 
& \leq \frac{(1 + \varepsilon)^{2/p}}{Y(M \times N, [g_r])}
\int_{M \times U}\bigl(a_{m+n}|\nabla f|_{g_r}^2 
+ {\bf s}_{h_r} f^2\bigr)d\mu_{g_r} \\ 
& \leq \frac{(1 + \varepsilon)^{2 + (2/p)}}{Y(M \times N, [g_r])}
\int_{M \times \mathbb{B}_{\rho}(0)} \left( a_{m+n}|\nabla f|_{\bar{g}}^2 
+ \left( {\bf s}_g + \frac{K_4}{r} \right)  f^2 \right) d\mu_{\bar{g}}, 
\end{align*} 

\noindent
where $K_4 > 0$ is a constant independent of $r$. 
Here, we also use the fact that 
$Y(M \times N, [g_r]) >0$ for any large $r > 0$. 

In order to prove (\ref{III}) it is important to note that for
$f \in C^{\infty}_c( \mathbb{B}_{\rho}(0))$, we obtain:
$$||f||^2_p
\leq \frac{(1 + \varepsilon)^{2 + (2/p)}}{Y_N(M \times N, g_r )}
\int_{M \times \mathbb{B}_{\rho}(0)} \left( a_{m+n}|\nabla f|_{\bar{g}}^2 
+ \left( {\bf s}_g + \frac{K_4}{r} \right) f^2 \right)d\mu_{\bar{g}}.$$ 

 From the positivity of ${\bf s}_g$, 
there exists $r_2 = r_2(\varepsilon, \underset{M}{\textrm{min}}{\bf s}_g ) > 0$ 
such that 
$${\bf s}_g + \frac{K_4}{r} \leq (1 + \varepsilon)
{\bf s}_g \quad \textrm{on}\ \ M$$ 
for any $r \geq r_2$. 
Hence, for any $r \geq \max\{r_1, r_2\}$ and 
$f \in C^{\infty}_c(M \times \mathbb{B}_{\rho}(0))$, 
we have 
$$||f||^2_p \leq 
\frac{(1 + \varepsilon)^{3 + (2/p)}}{Y(M \times N, [g_r])}
\int_{M \times \mathbb{B}_{\rho}(0)}\bigl(a_{m+n}|\nabla f|_{\bar{g}}^2 
+ {\bf s}_g f^2\bigr)d\mu_{\bar{g}} \ .$$ 
Therefore, 
$$Y(M \times N, [g_r]) \leq (1 + \varepsilon)^{3 + (2/p)}
Y(M \times \mathbb{B}_{\rho}(0), [\bar{g}]). $$ 
Letting $r \nearrow \infty$, we then obtain 
$$\limsup_{r \nearrow \infty}Y(M \times N, [g_r]) \leq 
(1 + \varepsilon)^{3 + (2/p)}Y(M \times \mathbb{B}_{\rho}(0), [\bar{g}]). $$ 
Letting also $\rho \nearrow \infty$ and $\varepsilon \searrow 0$, 
$$\limsup_{r \nearrow \infty}Y(M \times N, [g_r]) \leq Y_0.$$ 
And following the same steps we also prove (\ref{III}). 
This completes the proof of Theorem~\ref{3}. 
\end{proof}

We shall prove Theorem~\ref{4} 
by a series of lemmas. 
Throughout the rest of this section, 
we always assume the same conditions as in Theorem~\ref{4}, 
that is, $m, n \geq 2$ and ${\bf s}_g > 0$ on $M^m$. 
To simplify the notation, 
we set $\bar{g} := g + g_{{}_{\mathbb{E}}}$ 
on $M^m \times \re^n$. 

By the positivity of the scalar curvature ${\bf s}_g > 0$ of $g$ 
and the condition that $n \geq 2$, 
one can obtain the following. 

\begin{Lemma}\label{Weyl} 
$(M^m \times \re^n, \bar{g})$ is not locally conformally flat. 
Moreover, for any open set $U$ of $M^m \times \re^n$, 
the Weyl curvature $W_{\bar{g}}$ of $\bar{g}$ 
never vanishes identically on $U$. 
\end{Lemma} 

When $m + n \geq 6$, 
similarly to Aubin's argument in \cite{Aubin}, \cite[Theorem~B]{LP}, 
Lemma~\ref{Weyl} implies the existence of a family of nice test functions 
$\{\psi_{\varepsilon}\}_{\varepsilon > 0}$ 
with $Q_{\bar{g}}(\psi_{\varepsilon}) < {\bf Y}_{m + n}$ 
for sufficiently small $\varepsilon > 0$. 
Then we obtain the following 
(see the proof of Proposition~6.6 in \cite{AB} for details). 

\begin{Lemma}\label{h-dim} 
Assume that $m + n \geq 6$. 
Then, the assertion of Theorem~\ref{4} holds. 
\end{Lemma} 

Now we consider the remaining case when $m + n = 4, 5$ in Theorem~\ref{4}. 
Let $T^n_k := \re^n/(2^k\mathbb{Z})^n$ denote the quotient of $\re^n$ 
for $k = 0, 1, 2, \cdots$. 
Let $(\re^n, g_{{}_{\mathbb{E}}}) \rightarrow (T^n_k, h_k)$ 
also denote the natural infinite Riemannian covering 
and 
$$P_k : (M^m \times \re^n, \bar{g}) \rightarrow (M^m \times T^n_k, \bar{g}_k)$$ 
the induced Riemannian covering, where $\bar{g}_k := g + h_k$. 
Here note that each natural map 
$$(M^m \times T^n_{k +1}, \bar{g}_{k+1}) \rightarrow (M^m \times T^n_k, \bar{g}_k)$$ 
is also a finite Riemannian covering. 
Take any point $q \in M^m$ and fix it. 
Set $p = (q, 0) \in M^m \times \re^n$ 
and $p_k := P_k(p) \in M^m \times T^n_k$. 
Then, for each $k$, there exists a unique normalized positive Green's function $G_k$ 
(with pole at $p_k$) for the conformal Laplacian $L_{\bar{g}_k}$ 
on $M^m \times T^n_k$, that it, 
$G_k \in C^{\infty}_+\big((M^m \times T^n_k) - \{p_k\}\big)$ with 
$$L_{\bar{g}_k}G_k = c_{m + n}\ \delta_{p_k}\quad \textrm{on}\quad 
M^m \times T^n_k\qquad \textrm{and}\qquad 
\lim_{x \to p_k}\textrm{dist}(p_k, x)G_k(x) = 1.$$ 
Here, $c_{m+n} > 0$ and $\delta_{p_k}$ stand respectively for 
a specific univesal constant and the Dirac $\delta$-function at $p_k$. 
The condition ${\bf s}_g > 0$ implies that the first eigenvalue 
$\lambda(L_{\bar{g}_0}) > 0$ on $M^m \times T^n_0$. 
By the condition $\lambda(L_{\bar{g}_0}) > 0$, 
there exists a unique normalized minimal positive Green's function $G$ 
for $L_{\bar{g}}$ with pole at $p \in M^m \times \re^n$. 
Moreover, there exist positive constants $a, b, K$ with $a < b, K \geq 1$ such that 
$$K^{-1}e^{-b|z|} \leq G(x) \leq Ke^{-a|z|}\quad 
\textrm{for}\quad x = (y, z) \in M^m \times \{z \in \re^n\ \big|\ |z| \geq 1\}.$$ 
(See \cite{SY}, \cite[Section~6]{AB} for details.) 

Let us consider the family of the scalar-flat, asymptotically flat manifolds 
$$(X_k, \bar{g}_{k, AF}) := 
\big((M^m \times T^n_k ) - \{p_k\}, G_k^{\frac{4}{m+n-2}}\bar{g}_k\big)\quad 
\textrm{for}\quad k = 0, 1, 2, \cdots$$ 
and the one 
$$(X, \bar{g}_{AF}) := \big((M^m \times \re^n) - \{p\}, G^{\frac{4}{m+n-2}}\bar{g}\big)$$ 
with a singularity arising from the end of $M^m \times \re^n$. 
We denote the mass of $(X, \bar{g}_{AF})$ (resp. $(X_k, \bar{g}_{k, AF})$) 
by $\frak{m}(\bar{g}_{AF})$ (resp. $\frak{m}(\bar{g}_{k, AF})$). 
Then, a similar argument to the proof of \cite[Theorem~6.13]{AB} implies 
$$\lim_{k \to \infty}\frak{m}(\bar{g}_{k, AF}) = \frak{m}(\bar{g}_{AF}).$$ 
Hence, from the positive mass theorem \cite{SY-PM1, SY-PM2, Schoen}, 
we obtain the following. 

\begin{Lemma}\label{mass} 
Let $(X, \bar{g}_{AF})$ be the scalar-flat, asymptotically flat manifolds 
with a singularity and ${\rm dim}~X = 4, 5$, defined as abobe. 
Then, $\frak{m}(\bar{g}_{AF}) \geq 0.$ 
\end{Lemma} 

Now we can complete the proof of Theorem~\ref{4}. 

\begin{proof}[Proof of Theorem~\ref{4}] 
 From Lemma~\ref{h-dim}, we consider the case when $m + n = 4, 5$. 
We will prove that 
$$\frak{m}(\bar{g}_{AF}) > 0.$$ 
Then, by a similar argument to the proofs of \cite[Theorem~1]{Schoen0} 
and \cite[Chapter~5, Theorem~4.1]{SY-Book}, 
the positivity of the mass $\frak{m}(\bar{g}_{AF}) > 0$ 
implies the desired assertion 
$$Y(M^m \times \re^n, [\bar{g}]) < {\bf Y}_{m+n}.$$ 
Note also that, in order to prove $\frak{m}(\bar{g}_{AF}) > 0$, 
we may assume that $\bar{g}_{AF}$ is Ricci-flat on $X$ 
(cf.~\cite[Lemma~3]{Schoen0}, \cite[Proposition~6.14]{AB}). 

Now we suppose that $\frak{m}(\bar{g}_{AF}) = 0$. 
Choose a large constant $L_0 > 0$ and fix it. 
Set 
$$X_0 := \left( M^m \times \{z \in \re^n\ \big|\ |z| \leq L_0\} \right)
-\{ p \} \subset X.$$ 
Then, there exist {\it harmonic coordinates near infinity} 
$x = (x^1, \cdots, x^{m+n})$ on $(X_0, \bar{g}_{AF})$ 
\cite{Bartnik} (cf.~\cite[Lemma~6.17]{AB}). 
Namely, $(x^i)$ are smooth functions on $X_0$ satisfying 
$$ 
\Delta_{\bar{g}_{AF}}x^i = 0\quad \textrm{on}\quad X_0,\qquad 
\frac{\partial x^i}{\partial \nu} = 0\quad \textrm{on}\quad \partial X_0 
$$ 
and which give asymptotically flat coordinates near infinity 
of $(X_0, \bar{g}_{AF})$. 
Here, $\nu$ is the outward unit normal vector field normal to $\partial X_0$. 

We now apply the Bochner technique to completing the proof. 
The harmonicity of $(x^i)$ implies that $\{dx^i\}$ are harmonic $1$-forms 
on $(X_0, \bar{g}_{AF})$. From the Bochner formula for $1$-forms $\{dx^i\}$ 
combined with the conditions that 
$\frac{\partial x^i}{\partial \nu} = 0$ on $\partial X_0$ and 
$\textrm{Ric}_{\bar{g}_{AF}} = 0$ on $X$, we have that 
(cf.~\cite[Theorem~4.4]{Bartnik}, \cite[Proposition~10.2]{LP}) 
$$ 
\frak{m}(\bar{g}_{AF}) = 
\frac{1}{\textrm{Vol}\big(S^{m+n-1}(1)\big)}\sum_{i = 1}^{m+n}
\int_{X_0}|\nabla dx^i|^2d\mu_{\bar{g}_{AF}}. 
$$ 
Then, by applying the mass zero condition $\frak{m}(\bar{g}_{AF}) = 0$ 
in the above, 
we obtain that 
the $1$-forms $\{dx^i\}$ are parallel on $X_0$ with respect to $\bar{g}_{AF}$. 
Since the coframe $\{dx^i\}$ is orthonormal at infinity, 
then $\{dx^i\}$ is a parallel orthonormal coframe everywhere on $(X_0, \bar{g}_{AF})$. 
This implies that the map 
$x = (x^1, \cdots, x^{m+n}) : (X_0, \bar{g}_{AF}) \rightarrow (\re^{m+n}, g_{{}_{\mathbb{E}}})$ 
is a local isometry, 
and hence $\bar{g}$ is locally conformally flat on $X_0$. 
This gives a contradiction to Lemma~\ref{Weyl}. 
Therefore, $\frak{m}(\bar{g}_{AF}) > 0$. 
This completes the proof. 
\end{proof}

\section{Gagliardo-Nirenberg inequalities and Yamabe constants} 
In this section we estimate the behavior of arbitrary $N$--Yamabe constants 
in terms of the best constants in the Gagliardo-Nirenberg 
interpolation inequalities. 

\vspace{.06cm}

Let us define the 
{\it $(m,n)$--Gagliardo--Nirenberg functional} as 
$$L(f) = L_{m,n}(f) := 
\frac{{\| \nabla f \|}_2^{\frac{2n}{m+n}} {\| f \|}_2^{\frac{2m}{m+n}}}
{{\| f \|}_{p_{m+n}}^2}\qquad \textrm{for}\quad f\in W^{1,2}(\re^n) 
\quad \textrm{with}\quad f \not\equiv 0.$$ 

\begin{Remark}\label{inv} 
{\rm The map $L$ is invariant under two operations. First, if $c$ is
any non-zero constant then $L(cf) =L(f)$. Second, if $\lambda >0$ is
a constant and $f_{\lambda} (x)=f(\lambda x)$ then 
$${\| f_{\lambda} \|}_2^{\frac{2m}{m+n}}
={\lambda }^{\frac{-mn}{m+n}} {\| f \|}_2^{\frac{2m}{m+n}},$$ 
$$\, {\| f_{\lambda } \|}_{p_{m+n}}^2  = 
{\lambda}^{\frac{-2n}{p_{n+m}}}{\| f \|}_{p_{m+n}}^2$$ 
and 
$$\, \ \ {\| \nabla f_{\lambda} \|}_2^{\frac{2n}{m+n}} = 
{\lambda}^{\frac{(2-n)n}{m+n}} {\| \nabla f \|}_2^{\frac{2n}{m+n}}.$$ 
Therefore, $L(f_{\lambda} )=L(f)$.} 
\end{Remark} 

\noindent 
Let us recall the following definition from the introduction: 

\begin{Definition}\label{sob1} 
The $(m,n)$--Gagliardo--Nirenberg constant is given by
$${\sigma}_{m,n} := 
{\left( \inf_{f\in W^{1,2}({\re}^n) \atop f \not\equiv 0} L_{m,n}(f)\right)}^{-1}.$$ 
\end{Definition} 

\begin{Remark}\label{num} 
{\rm The constant $\sigma_{m,n}$ has 
already been studied in the literature. It is the best constant 
for the classical interpolation inequality due to Gagliardo and 
Nirenberg that says that $L_{m,n}$ is bounded away from zero 
(cf \cite{g1}, \cite{g2} and \cite{n}). In \cite{Weinstein}, 
it is shown that $\sigma_{m,n}$ is closely related to the global
existence of the nonlinear Schr\"odinger equation. 
The author showed also that 
$\sigma_{m,n}$ is always attained by a positive function 
$\psi \in W^{1,2}(\re^n) \cap  C^\infty(\re^n)$, called a 
{\it ground state solution}, that should then satisfy the 
corresponding Euler--Lagrange equation 
\begin{equation}\label{gss} 
\Delta u-u+u\sp q=0, \ \ \ \ q={\frac{m+n+2}{m+n-2}}
\end{equation} 
(the Euler-Lagrange equation for $L$ of course involves coefficients
depending on ${\|f \|}_2$, ${\| f \|}_p$
and $ {\| \nabla f \|}_2$: one can use the previous remark to normalize 
the equation as above). 
For any function $f\in W^{1,2}(\re^n) \cap  C^\infty(\re^n)$, 
one can consider the spherical symmetrization $f^*$ of $f$, that is, 
the radial function $f^*$ satisfying 
$\textrm{Vol}_{g_{{}_{\mathbb{E}}}}\{ f^* > t \} 
= \textrm{Vol}_{g_{{}_{\mathbb{E}}}}\{ f > t \}$ 
for any positive $t$. 
It is a classical result that $L(f^* ) < L(f)$ 
if $f^* \not\equiv f$. 
It follows that $\psi$ should be radial and decreasing. 
Finally, in \cite{k} it is proved the uniqueness of the positive 
radial solution of (\ref{gss}) under the assumption that it  
vanishes at infinity. The key point to our purposes is 
that these facts give a simple procedure to compute numerically 
all Gagliardo--Nirenberg constants. We will continue this discussion 
on Section 6.} 
\end{Remark} 

\vspace{0.1cm} 

We now relate the Gagliardo--Nirenberg constants 
with the $N$--Yamabe constant of limiting products. 

{\it Proof of Theorem \ref{1}}. 
First fix a function 
$f \in C^{\infty}_c( \mathbb{R}^n )$ 
and consider the family of functions $f_{\lambda}$ as in Remark \ref{inv}. 
Let us consider the map
$$\lambda \mapsto F(\lambda ) := 
\frac{\int_{\re^n} \big(a_{m+n} |\nabla f_{\lambda} |^2_{g_{{}_{\mathbb{E}}}}
+{\bf s}_g f_{\lambda}^2\big)  \  d\mu_{g_{{}_{\mathbb{E}}}} }
{ ||f_{\lambda} ||^2_{p_{m+n}} } \ .$$ 
\noindent 
If we set 
$$A := 
\frac{\int_{\re^n} a_{m+n} |\nabla f |^2_{g_{{}_{\mathbb{E}}}}d\mu_{g_{{}_{\mathbb{E}}}}}
{ ||f ||^2_{p_{m+n}} }$$ 
\noindent 
and
$$B := 
\frac{\int_{\re^n} {\bf s}_g f^2  \  d\mu_{g_{{}_{\mathbb{E}}}} }{ ||f ||^2_{p_{m+n}} } \ ,$$ 
\noindent 
then Remark \ref{inv} implies that 
$F(\lambda )={\lambda}^{\frac{2m}{m+n}} A + \lambda^{\frac{-2n}{m+n}} B$. 
We see that this map achieves its minimum at 
$\lambda_0 =  \sqrt{\frac{nB}{mA}}$ and that this minimum is 
$$F(\lambda_0 ) = 
A^{\frac{n}{m+n}}B^{\frac{m}{m+n}}m^{\frac{-m}{m+n}}n^{\frac{-n}{m+n}} (m+n) 
= {\bf s}_g^{\frac{m}{m+n}} C(m,n) L(f).$$ 
\noindent 
Recall that $L(f) = L(f_{\lambda })$. Therefore, 
$$\lim_{r \nearrow \infty} Y_{N^n}(M^m \times N^n, g + r h ) 
= Y_{\re^n}(M^m \times \re^n, g + g_{{}_{\mathbb{E}}})$$ 
$$ = \inf_{f \in C_c^{\infty}(\re^n ) \atop f \not\equiv 0}\ 
\frac{\int_{\re^n} \big(a_{m+n} |\nabla f |^2_{g_{{}_{\mathbb{E}}}}
+{\bf s}_g f^2\big)  \  d\mu_{g_{{}_{\mathbb{E}}}} }
{ ||f ||^2_{p_{m+n}} } $$ 
$$ = \inf_{ f \in C^{\infty}_c( \mathbb{R}^n ) \atop f \not\equiv 0}
{\bf s}_g^{\frac{m}{m+n}}C(m,n)L(f) 
= \frac{{\bf s}_g^{\frac{m}{m+n}} C(m,n) }{\sigma_{m,n}} \ .
$$ 

\qed

\section {Numerical computations} 
We describe now a procedure to determine all Gagliardo--Nirenberg 
constants $\sigma_{m,n}$ numerically, and how we obtained our table 
in the introduction.

\medskip 

Consider the solution $h = h_\alpha:\re_{\geq 0} \to \re$ of the initial value problem 

\begin{equation}\label{ode} 
h''(t) + \frac{n-1}{t}h'(t) - h(t) + h(t)^{\frac{m+n+2}{m+n-2}} = 0,
\ \ \ h(0) = \alpha > 0,\ h'(0) = 0. 
\end{equation} 
This equation corresponds to the critical points of the
$(m,n)$--Gagliardo--Nirenberg functional $L_{m,n}$ 
(see Remarks~\ref{inv}, \ref{num}), 
when restricted to radial functions 
$u(x) = h(|x|)$, which is enough to consider because of symmetrization. 
By the existence result in \cite{Weinstein} and its uniqueness finally proved in \cite{k}, 
there exists only one value 
$\alpha = \alpha_0 = \alpha_0(m,n)$ that gives a positive solution 
$h_{\alpha_0}$ that vanishes at infinity, called the 
{\it ground state}. 
To find $\alpha_0$ numerically, we use \cite{k} where it is 
shown that, for values $\alpha > \alpha_0$, the solution $h_\alpha$ 
vanishes exactly once and then oscillates about $-1$ (Figure~1), 
while, for values $\alpha < \alpha_0$, it is positive and oscillatory 
about the value $1$ (Figure~2). 

\epsfysize=3 cm

\begin{figure}[tbp]
\center{
\epsfbox{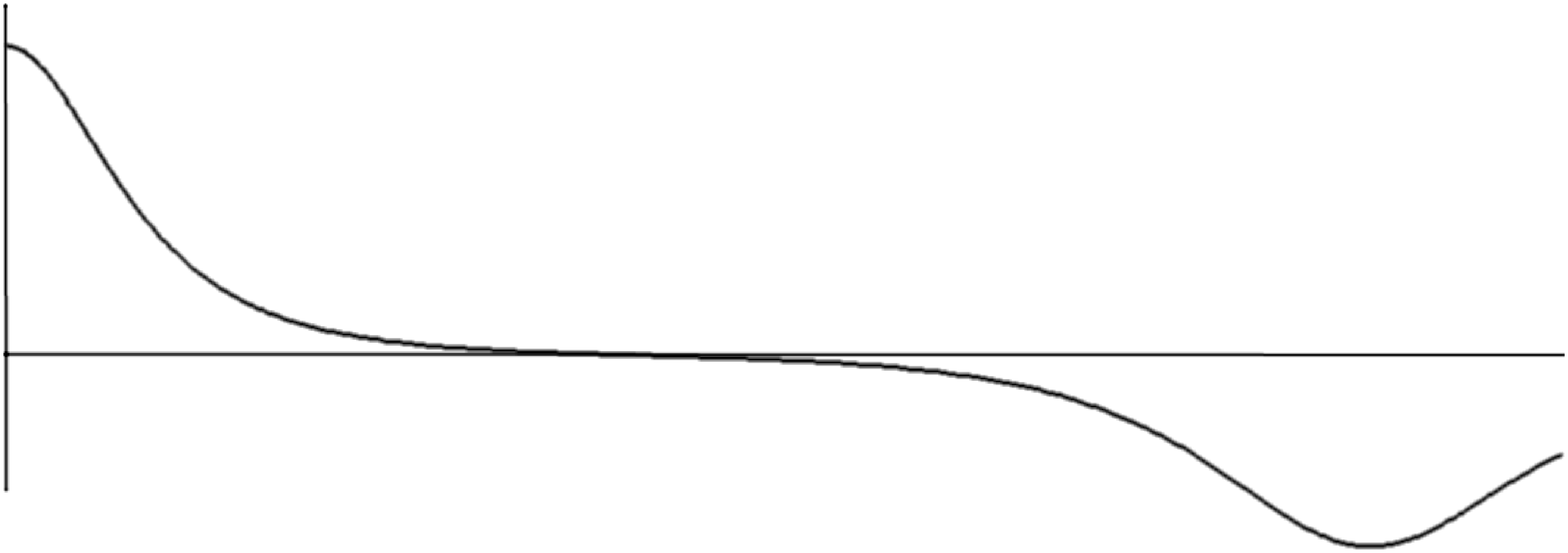}
\vskip -0.3cm
{\tiny Figure 1: {$h_\alpha$ for $\alpha>\alpha_0$.}}
}
\end{figure}

\epsfysize=3 cm

\begin{figure}[tbp]
\center{
\epsfbox{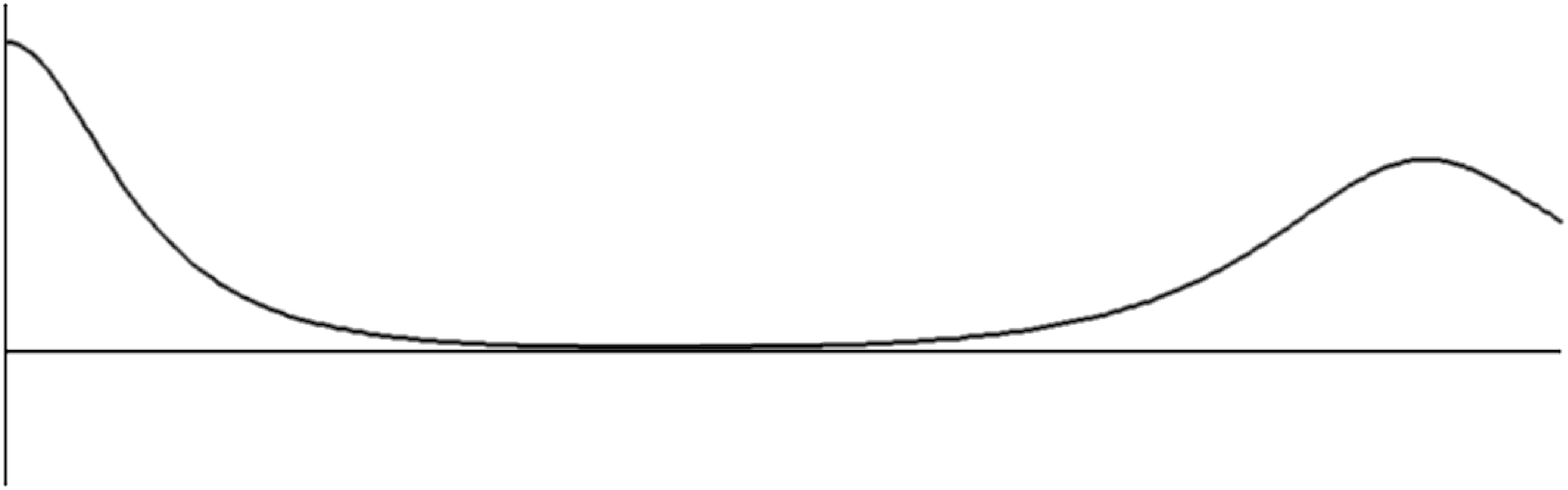}
\vskip -0.3cm
{\tiny Figure 2: {$h_\alpha$ for $\alpha<\alpha_0$.}}
}
\end{figure} 
A key point is the uniqueness of the solution of (\ref{ode}) 
when the initial value condition $h(0) = \alpha$ is replaced by 
the boundary condition $h(t_0) = 0$, for $t_0\in (0,+\infty]$, and 
the fact that this solution has $h'<0$ in $(0,t_0]$; see Lemmas 9, 11 
and Theorem in \cite{k}. 
Finally, by \cite{Weinstein} we have that 
$\sigma_{m,n}^{-1}= L_{m,n}(h_{\alpha_0})$. 
For example, we can then 
compute $\alpha_0 = \alpha_0(2,2) \thickapprox 2.2062$ 
for the ground state initial value, and hence 
$\sigma_{2,2}^{-1}= L_{2,2}(h_{\alpha_0}) \thickapprox 2.41877$. 
In fact, Figure 1 corresponds to $m=n=2$ and $\alpha=2.208$, while 
Figure~2 to $\alpha=2.205$. 

Of course, if one wants to avoid the numerical computation, 
one could give estimates for $\sigma_{m,n}$ by carefully choosing functions. 
For instance we can show that 
$\sigma_{2,2}^{-1} < 2.427458 <\sqrt{2\pi}$. 
by considering the function $h:\re_{\geq 0} \to \re$ 
that linearly interpolates the following data: 
$h(0)    = 1,$ 
$h(0.1)  = 0.9904132,$ 
$h(0.2)  = 0.9626,$ 
$h(0.3)  = 0.91917,$ 
$h(0.7)  = 0.66607,$ 
$h(0.9)  = 0.5378,$ 
$h(1.15) = 0.4023,$ 
$h(1.3)  = 0.34$, 
$h(1.5)  = 0.2634,$ 
$h(1.85) = 0.17201,$ 
$h(2.2)  = 0.11288,$ 
$h(2.6)  = 0.07031,$ 
$h(3)    = 0.04416,$ 
$h(3.5)  = 0.02493,$ 
$h(3.9)  = 0.016,$ 
$h(4.3)  = 0.01016,$ 
$h(5)    = 0.0047,$ 
$h(6)    = 0.00158,$ 
$h(7)    = 0.00054,$ 
$h(8)    = 0.00019,$ 
$h(9)    = 0.00006,$ and 
$h(t)    = 0$, for $t\geq 10$. 
Now, define the radial function $f:\re^2\to\re$ given by 
$f(x)=h(|x|)$. A straightforward computation gives that 
$\sigma_{2,2}^{-1} \leq L_{2,2}(f)<2.427458$. 
And then 
$$Y_{S^2}^{\infty}(S^2 \times S^2) < 
8\sqrt{3\pi} \times 2.427458 < 8\sqrt{6}~\pi = {\bf Y}_4.$$

\vspace{0.5cm}

\end{document}